\documentclass[12pt]{amsart}
\usepackage{amssymb,amsmath,latexsym}

\newcommand\comment[1]{}				
\renewcommand\comment[1]{\emph{[#1]}}			

\newcommand\cA{{\mathcal A}}
\newcommand\cM{{\mathcal M}}
\newcommand\cP{{\mathcal P}}
\begin{document}

\begin{center}

\Large{\bf A Shorter, Simpler, Stronger Proof of the \\ 
           Meshalkin--Hochberg--Hirsch Bounds \\ 
           on Componentwise Antichains\footnote{Appeared in \emph{Journal of Combinatorial Theory Series A} {\bf 100} (2002), 196--199.}}
\normalsize

{\sc Matthias Beck and Thomas Zaslavsky\footnote{Research supported by National Science Foundation grant DMS-70729.}}

{\sc State University of New York at Binghamton \\
     Binghamton, NY, U.S.A.\ 13902-6000}

{\tt matthias@math.binghamton.edu \\
     zaslav@math.binghamton.edu }

\end{center}

\bigskip 

{\it Abstract:} 
{Meshalkin's theorem states that a class of ordered $p$-partitions of an $n$-set 
has at most $\max \binom{n}{a_1,\hdots, a_p}$ members if for each $k$ the 
$k^{\textrm{th}}$ parts form an antichain.  We give a new proof of 
this and the corresponding LYM inequality due to Hochberg and Hirsch, which is 
simpler and more general than previous proofs.  It extends to a common generalization of Meshalkin's theorem and Erd\H{o}s's theorem about $r$-chain-free set families.}
\normalsize



\emph{Running head}: Meshalkin--Hochberg--Hirsch Bounds 

\emph{Address for editorial correspondence}:\\
    Thomas Zaslavsky\\
    Department of Mathematical Sciences\\
    State University of New York\\
    Binghamton, NY 13902-6000 \\
    U.S.A.

\vfill\pagebreak

An \emph{antichain of sets} is a class of sets of which none contains another. 
Sperner \cite{sperner} bounded the size of an antichain $\cA$ of subsets 
of an $n$-element set $S$: 
$$ 
|\cA| \leq {n\choose {\lfloor n/2 \rfloor }} \ , 
$$ 
with equality if $\cA = \cP_{\lfloor n/2 \rfloor} (S)$ or $\cP_{\lceil n/2 \rceil}(S)$, 
where by $\cP_k(S)$ we mean the class of $k$-element subsets of $S$. 
Subsequently, Lubell \cite{lubell}, Yamamoto \cite{yamamoto}, and 
Meshalkin \cite{meshalkin} independently obtained a stronger result 
(of which Bollob\'as independently proved a generalization \cite{bollobas}): 
any antichain $\cA$ satisfies 
$$ 
\sum_{A\in \cA}\ \frac {1}{\binom n {|A|}} \leq 1 \ , 
$$ 
and equality holds if $\cA = \cP_{\lfloor n/2 \rfloor} (S) $ or $\cP_{\lceil n/2 \rceil}(S)$. 


We give a very short proof of a considerable generalization of these results. 
A class $\cA$ of subsets of $S$ is \emph{$r$-chain-free} if 
$\cA$ contains no chain of length $r$.\footnote{An $r$-chain-free family has been called an ``$r$-family'' or 
``$k$-family'', depending on the name of the forbidden length, but we think it is time for a distinctive name.} 
(A \emph{chain} is a class of mutually comparable sets, that is, 
$T\subset T' \subset \cdots \subset T^{(l)}$. Its \emph{length} is $l$.) 
A \emph{weak composition of $S$ into $p$ parts} 
is an ordered $p$-tuple $A = (A_1,\ldots,A_p)$ such that the $A_k$ are 
pairwise disjoint subsets of $S$ and their union is $S$.  We call $A_k$ the 
\emph{$k^{\textrm{th}}$ part} of $A$. A part $A_k$ may be void 
(hence the word ``weak'').  
If $\cM = \{ A^1, \hdots, A^m\}$ is a class of weak compositions of $S$ into $p$ parts, 
we write $\cM_k = \{ A^i_k\}^m_{i=1}$ for the class of distinct $k^{\textrm{th}}$ parts of members of $\cM$. 
A multinomial coefficient of the form $\binom{n}{a_1, \hdots, a_p}$ is called a 
\emph{$p$-multinomial coefficient for $n$}. Our result is: 

\noindent {\bf Theorem.} 
{\it 
Let $n\geq 0$, $r\geq 1$, $p\geq 2$, and let $S$ be 
an $n$-element set.  Suppose $\cM$ is a class of weak compositions of $S$ into $p$ 
parts such that, for each $k<p$, $\cM_k$ is $r$-chain-free. Then 

\noindent {\rm (a)} $ \displaystyle \sum_{A\in \cM} \dfrac{1}{\binom{n}{|A_1|,\hdots,|A_p|}} \leq r^{p-1} $, and 

\noindent {\rm (b)} $|\cM|$ is bounded by the sum of the $r^{p-1}$ largest $p$-multinomial coefficients for $n$.

} 

The number of $p$-multinomial coefficients for $n$ is $\binom{n+p-1}{p-1}$; 
if $r^{p-1}$ exceeds this we extend the sequence of coefficients with zeros.

Our theorem is a common generalization of results of Meshalkin and Erd\H{o}s.  
The case $r=1$ of part (b) (with the added assumption that every $\cM_k$ is an 
antichain) is the relatively neglected theorem of Meshalkin \cite{meshalkin}; 
later Hochberg and Hirsch \cite{hh} found (a) for this case, which implies (b). 
Our extension to $r>1$ is inspired by the case $p=2$, which is equivalent to 
Erd\H{o}s's theorem \cite{erdos} that 
for an $r$-chain-free family $\cA$ of subsets of $S$,
$|\cA |$ is bounded by the sum of the $r$ largest binomial coefficients ${n\choose k}$, $0\leq k \leq n$, and its LYM companion due to Rota and Harper \cite{rota}. 
We need the latter for our theorem; we sketch its proof for completeness' sake.


\noindent {\bf Lemma} \cite[p.\ 198, ($^*$)]{rota}. 
{\it 
For an $r$-chain-free family $\cA$ of subsets of $S$,
$$ 
\sum_{A\in \cA} \frac {1}{\binom n {|A|}} \leq r \ . 
$$ 
} 

\begin{proof} 
Each of the $n!$ maximal chains in $\cP(S)$ contains at most $r$ members of $\cA$. 
On the other hand, there are $|A|! (n-|A|)!$ maximal chains containing $A \in \cP(S)$. 
Now count: 
$$ 
\sum_{A\in \cA} |A|! (n-|A|)! \leq r n! \ . 
$$ 
The lemma follows.
\end{proof}

\begin{proof}[Proof of the theorem] 
Our proof of our whole theorem is simpler than the original proofs of the case $r=1$. 
The proof of (a) here (which is different from the more complicated although equally short proof by Hochberg and Hirsch) is inspired by the beginning of Meshalkin's proof of (b) for $r=1$. 
We proceed, as did Meshalkin, by induction on $p$.  The case $p=2$ is 
equivalent to the lemma because 
if $(A_1, A_2)$ is a weak composition of $S$, then $A_2 = S\setminus A_1$.  Suppose 
then that $p>2$ and (a) is true for $p-1$.  Let 
$\cM(F) = \{ (A_2, \hdots, A_p): (F, A_2, \hdots, A_p) \in \cM \}$. 
Since $\cM(F)_k \subseteq \cM_{k+1}$, $\cM(F)_k$ is $r$-chain-free for $k<p-1$. Thus, 
\begin{align*}
\sum_{A\in \cM} \frac {1}{\binom{n}{|A_1|, \hdots, |A_p|}}
&= \sum_{A\in \cM} \frac{1}{\binom{n}{|A_1|}} 
    \frac{1}{\binom{n-|A_1|}{|A_2|, \hdots, |A_p|}}  \\
&= \sum_{F\in \cM_1} \frac{1}{\binom{n}{|F|}} \sum_{A' \in \cM (F)} 
    \frac{1}{\binom{n-|F|}{|A_2|, \hdots, |A_p|}}
\intertext{where $A' = (A_2,\hdots,A_p)$,}
&\leq \sum_{F\in \cM_1} \frac{1}{\binom{n}{|F|}} r^{p-2}
\intertext{by the induction hypothesis,}
&\leq r \cdot r^{p-2}
\end{align*}
by the lemma. This proves (a).

To deduce (b), write the $p$-multinomial coefficients for $n$ in any weakly decreasing order as $M_1, M_2, \hdots$, 
extended by $0$'s as necessary to a sequence of length $r^{p-1}$.  
In the left-hand side of (a), replace each of the $M_1$ terms with largest denominators by $1/M_1$.  Their
sum is now $1$.  Amongst the remaining terms all denominators are at most $M_2$; replace the $M_2$ of them 
with the largest denominators by $1/M_2$. Now their sum is $1$.  Continue in this fashion.  The number of 
terms could be less than $T = M_1+ \cdots + M_{r^{p-1}}$; in that case, $|\cM|<T$.  Otherwise, after $r^{p-1}$ 
steps we have replaced $T$ terms and have on the left side of (a) a sum equal to $r^{p-1}$ plus any further 
terms. As the total is no more than $r^{p-1}$, there cannot be more than $T$ terms.  Thus we have proved (b).
\end{proof}

Our proof is naturally general: that of (a) would be no shorter even if restricted to $r=1$ (but the deduction of (b) would become trivial).
What is more, it is applicable to projective geometries \cite{bz}. 
Furthermore, our proof, even restricted to $r=1$, is simpler than the original proofs by Meshalkin and Hochberg--Hirsch.

The upper bounds in the theorem can be attained only in limited circumstances. 
When $r=1$, the maxima are attained if for each $k$, $\cM_k = \cP_{\lfloor n/p \rfloor } (S)$ or $\cP_{\lceil n/p \rceil}(S)$ \cite{meshalkin}. 
When $p=2$, the upper bounds are attained if $\cM_1$ is the union of the $r$ largest classes $\cP_m(S)$ \cite{erdos}.
When $r>1$ and $p>2$, the upper bounds are only sometimes attainable, but proving this is complicated.



\bibliographystyle{alpha}

\end{document}